\documentclass{amsart}

\usepackage{amsmath, amsthm, amssymb, latexsym,url}
\usepackage{hyperref}

\title{On the binary digits of $\sqrt{2}$}
\author{Joseph Vandehey}

\newtheorem{thm}{Theorem}[section]

\newtheorem{rem}[thm]{Remark}

\begin{document} 

\maketitle

\begin{abstract}
We show that the number of $1$'s in the first $N$ digits of the binary expansion of $\sqrt{2}$ is at least $\sqrt{2N}(1+o(1))$ and show that this bound can be improved to around $2\sqrt{N}/\sqrt{2\sqrt{2}-1}$ infinitely often.
\end{abstract}

\section{Introduction}

It has long been a folklore conjecture that $\sqrt{2}$ is a normal number in any given base---that is, the limiting frequency of any string of digits in the base-$b$ expansion of $\sqrt{2}$ is the same as any other string of the same length. If $\sqrt{2}$ were normal in base $2$, then we would expect that
\begin{equation}\label{eq:weak}
\operatorname{nz}(N)=\frac{N}{2}(1+o(1))
\end{equation}
where $\operatorname{nz}(N)$ is the number of non-zero digits in the first $N$ binary digits after the decimal point of $\sqrt{2}$.

However, we cannot yet prove a result of the strength of \eqref{eq:weak}. Indeed, the strongest result we have available is
\[
\operatorname{nz}(N)\ge N^{1/2}(1+o(1)),
\]
due to Bailey, Borwein, Crandall, and Pomerance \cite{BBCP}. Although their result was significantly more general (looking at arbitrary algebraic irrationals and arbitrary bases), and although several papers have been written improving the bounds in these general cases \cite{AF,Kaneko,Riv}, the bound in this ``simplest" case of $\sqrt{2}$ in base $2$ has been unchanged. A survey article by Kaneko \cite{Kaneko2} has many related results.

Our main theorem of this paper is the following.

\begin{thm}\label{thm:main}
As $N$ tends to infinity, we have
\[
\operatorname{nz}(N) \ge \sqrt{2}N^{1/2} (1+o(1))
\]
\end{thm}

Our proof will largely follow the methods of Bailey, et al, but we also owe thanks to Bugeaud's exposition of their proof in his book \cite{Bugeaud}.

We can offer occasional improvements to the above theorem.

\begin{thm}\label{thm:2}
For any $\epsilon>0$ and infinitely many $N$, we have
\[
\operatorname{nz}(N)\ge \left(\frac{2}{\sqrt{2\sqrt{2}-1}}-\epsilon\right) N^{1/2}.
\]
\end{thm}

We note that $\sqrt{2}\approx 1.41421$, while $2/\sqrt{2\sqrt{2}-1}\approx 1.47908$. And the constant here could be improved further, possibly up to $\sqrt{8/\pi}\approx 1.59577$, but in the interest of expediency, we only include the simpler proof.

We conclude the paper with two ideas, which, while interesting, do not improve the results given above. In section \ref{sec:oddeven}, we look at the nonzero digits $a_i$ along odd and even indices $i$. In section \ref{sec:3root2}, we compare the digits of $\sqrt{2}$ and $3\sqrt{2}$ and see they can't simultaneously both be too small. 

We use standard asymptotic notation in this paper. By $f(x)=O(g(x))$ we mean that there exists a constant $C>0$ such that $|f(x)|\le C\cdot |g(x)|$. By $f(x)=o(g(x))$ we mean that as $x\to \infty$, $f(x)/g(x)$ tends to $0$.

\section{The $T$ and $r$ functions}

Again, we borrow the notation from the Bailey, et al, paper, with minor changes.

We consider $a_i\in\{0,1\}$ to be the binary digits of $\sqrt{2}$, so that 
\[
\sqrt{2}= \sum_i \frac{a_i}{2^i}
\]

We let 
\[
r(n) = \#\{(i,j):i+j=n,\ a_i=1, \ a_j=1\}.
\]
In particular, this definition means that
\[
2=\left( \sum_i \frac{a_i}{2^i} \right)^2 = \sum_n \frac{r(n)}{2^n}.
\]

We now also define
\begin{equation}\label{eq:firstTeq}
T(R)= \sum_{m\ge 1} \frac{r(m+R)}{2^m}.
\end{equation}
as the ``tail component." We may think of this in the following way: when squaring $\sqrt{2}$ to obtain $2$ digit by digit, the $T(R)$ function measures the contribution from digits at the $R+1$st place, including any carries that occur.

Now it is clear that $r$ is always a non-negative function, and must be positive infinitely often. Hence, $T(R)$ is always positive.

Moreover, it is a trivial consequence of \eqref{eq:firstTeq} that 
\begin{equation}\label{eq:secondTeq}
2T(R-1) = T(R)+r(R).
\end{equation}
Thus, in particular, since $T(-1)=1/2$ is a half-integer and $r(R)$ is always integer-valued, $T(R)$ is always an integer for all non-negative values of $R$.

The key idea of the original proof of Bailey, et al, is to have a lower bound and an upper bound on
$\sum_{R=1}^N T(R)$, and the upper bound relies on having a certain
size of $\operatorname{nz}(N)$. 

\section{A simple lower bound on $T(R)$}\label{sec:lower}

Consider the facts that $r$ is always non-negative, $T$ is always positive, and
\eqref{eq:secondTeq}. Combined, these tell us that
\[
T(R-1)> \frac{r(R)}{2},
\]
and thus
\begin{equation}\label{eq:T(R-1)bound}
T(R-1) \ge \begin{cases}
\frac{r(R)}{2} +1& \text{if }r(R) \text{ is even}\\
\frac{r(R)}{2} +\frac{1}{2}& \text{if }r(R) \text{ is odd}\\
\end{cases}
\end{equation}
We expect that odd values of $r(R)$ to be a (reasonably) rare
event, since, by symmetry, the only way for $r(R)$ to be odd is if
$R=2m$ and $a_m=1$. This happens at most $\operatorname{nz}(N)$ times up to
$N$.

Thus, we have
\begin{equation}\label{eq:basiclower}
\sum_{R=0}^{N-1} T(R) \ge \sum_{R=1}^{N} \left(  \frac{r(R)}{2}
  +1\right) + O\left( \operatorname{nz}(N)\right) = \frac{1}{2}
\sum_{R=1}^N r(R) + N + O\left( \operatorname{nz}(N)\right).
\end{equation}

\section{A simple upper bound}

Suppose that $N$ is sufficiently large and let $K=K(N)$ be a function dependent on $N$ to be defined shortly.

By \eqref{eq:firstTeq} we have
\begin{align*}
\sum_{R=0}^{N-K} T(R) &=\sum_{R=0}^{N-K} \sum_{m \ge 1} \frac{r(m+R)}{2^m}\\
&\le  \sum_{R=1}^\infty r(R) \cdot \sum_{m=\max\{R-N+K,1\}}^\infty \frac{1}{2^m}\\
&\le \sum_{R=1}^N r(R) + \sum_{R=N+1}^\infty \frac{r(R)}{2^{R-N+K}}.
\end{align*}

Now we make use of the fact that by definition, $r(n)\le n+1$. So
\[
\sum_{R=0}^{N-K} T(R)\le  \sum_{R=1}^N r(R) + \frac{N+3}{2^K}.
\]
So if we choose $K=\lfloor \log_2 N\rfloor$, then
\begin{equation}\label{eq:firstupper}
\sum_{R=0}^{N-K} T(R) \le \sum_{R=1}^N r(R) +O(1).
\end{equation}

\section{Proving $\operatorname{nz}(N)\ge  N^{1/2}(1+o(1))$}

Now, we can reprove the bound of Bailey, et al. It is clear
that \begin{equation}\label{eq:basicupper}
\sum_{R=1}^N r(R) \le (\operatorname{nz}(N))^2,
\end{equation} since the latter can be interpreted as $\#\{(i,j):i\le N,j\le N,a_i=a_j=1\}$, whereas the former can be interpreted as $\#\{(i,j):1\le i+j\le N, a_i=a_j=1\}$. Combining
this with \eqref{eq:basiclower} and \eqref{eq:firstupper}, we get
\[
N-K+O\left( \operatorname{nz}(N-K)\right) \le (\operatorname{nz}(N))^2.
\]
Noting that $\operatorname{nz}(N)$ is a non-decreasing function that tends to infinity with $N$ and that $K=o(N)$, we have
\[
N^{1/2}(1+o(1)) \le \operatorname{nz}(N).
\]

\section{The proof of Theorem \ref{thm:main}}

First, note that $T(R)$ is a positive
integer and expressible as $2T(R-1)-r(R)$. Thus $T(R)$ is only odd if
$r(R)$ is odd, and $r(R)$ is only odd if $R=2m$ for some $m$ with $a_m=1$. So $T(R)\ge 2$ unless $R$ is twice the index $i$ of some non-zero digit $a_i$ of which there are at most $\operatorname{nz}(N)$ many up to $N$, and thus we get 
\begin{equation}\label{eq:refinedlower}
\sum_{R=1}^N T(R) \ge 2N + O\left( \operatorname{nz}(N)\right).
\end{equation}
Again combining this with \ref{eq:firstupper} and \ref{eq:basicupper}
gives the bound
\[
\sqrt{2}N^{1/2} (1+o(1))\le \operatorname{nz}(N).
\]

\section{The proof of Theorem \ref{thm:2}}

Our upper bound \eqref{eq:basicupper} is suboptimal because it implicitly assumes that every pair $(i,j)$ with $a_i=a_j=1$ and $i,j\le N$ also satisfies $i+j \le N$. In fact, if $i,j\in (N/2,N]$ then they cannot possibly contribute to the sum of $r(R)$'s with $R\le N$. Thus, we could consider the improved upper bound
\begin{equation}\label{eq:refinedupper}
\sum_{R=1}^N r(R) \le (\operatorname{nz}(N/2))^2 + 2\operatorname{nz}(N/2)\left(\operatorname{nz}(N)-\operatorname{nz}(N/2)\right)
\end{equation}
This bound comes about because each sum $i+j\le N$ must come about because either both $i,j$ are in $[0,N/2]$ or one of them is in $[0,N/2]$ and the other is in $(N/2,N]$.

Suppose that Theorem \ref{thm:2} is not true, so there exists some $\epsilon>0$ such that for all sufficiently large $N$,
\[
\operatorname{nz}(N)< \left( \sqrt{2/(2\sqrt{2}-1)}-\epsilon\right) N^{1/2}.
\]
Let $\lambda = 2/(2\sqrt{2}-1)-\epsilon$.

Then, by a standard argument, we can find functions $g_1(N)$, $g_2(N)$ that are both $o(N)$, positive constants $\lambda_1,\lambda_2$ that are both at most $\lambda$, and an increasing sequence $\{N_i\}_{i=1}^\infty$ of positive integers such that
\begin{align*}
\operatorname{nz}(N/2)&= \lambda_1 \sqrt{N}+g_1(N)\\
\operatorname{nz}(N) &= \lambda_2 \sqrt{2N}+g_2(N),
\end{align*}
for $N$ belonging to the sequence of $N_i$'s.

By using \eqref{eq:firstupper}, \eqref{eq:refinedlower}, and \eqref{eq:refinedupper}, we see that for each $N$ belonging to the sequence of $N_i$'s, we have
\begin{align*}
2(N-K(N) )+O(\operatorname{nz}(N-K(N)))&\le (\lambda_1^2 N+ 2 \lambda_1(\lambda_2\sqrt{2}-\lambda_1) N)(1+o(1))\\ &= \lambda_1 (2\sqrt{2} \cdot \lambda_2 - \lambda_1)N(1+o(1)).
\end{align*}
We may assume that $\operatorname{nz}(N-K(N))=o(N)$ as otherwise we obtain a much stronger result than Theorem \ref{thm:2}. Therefore, we may simplify this inequality to
\[
2\le \lambda_1 (2\sqrt{2} \cdot \lambda_2 - \lambda_1),
\]
since the contributions of $(1+o(1))$ are no longer relevant without $N$ in the equation.

Since $\lambda_2\le \lambda$, we have that
\[
2\le \lambda_1 (2\sqrt{2}\cdot \lambda-\lambda_1).
\]
Moreover, by taking the derivative with respect to $\lambda_1$ on the right-hand side, we see that this is an increasing function in $\lambda_1$ for $\lambda_1\le \sqrt{2}\lambda$, and since $\lambda_1\le \lambda$, we have
\[
2\le (2\sqrt{2}-1)\lambda^2.
\]
However, this is a clear contradiction to the definition of $\lambda$.

\begin{rem}

This method can be extended considerably. We used only two intervals $[0,N/2]$ and $(N/2,N]$. By breaking into $m$ intervals and letting $m$ tend to infinity, it appears that we can show that for any $\epsilon>0$ there exist infinitely many integers $N$ such that
\[
\operatorname{nz}(N) \ge \left(\sqrt{\frac{8}{\pi}}-\epsilon\right)N^{1/2}.
\]
 The appearance of the $\sqrt{\pi}$ is due to the resulting sums coming closer and closer to an integral that resembles $\sqrt{1/x(1-x)}$.  However, given the increased difficulty of the proof and the negligible improvement it offers, we do not write it here.
\end{rem}

\section{Odd and even indices}\label{sec:oddeven}

There's no necessary reason why we have to look at $\sum_{R=1}^N T(R)$ on its own. We could weight this sum or look at sums along certain sequences. For example, we could look along arithmetic progressions, such as all even numbers. If we did that, we  would get a result that looks like
\[
\sum_{R=0}^{N-K} T(2R) < \sum_{R=1}^N \left( r(2R+1)+ \frac{1}{2} r(2R)\right) + o(1).
\]

Now suppose we let $\operatorname{nz}_0(N)$ denote the number of $i\le 2N+1$ such that $i$ is even and $a_i=1$, and we let $\operatorname{nz}_1(N)$ denote the number of $i\le 2N+1$ such that $i$ is odd and $a_i=1$. Then it is clear that
\[
\sum_{R=0}^N r(2R) \le \operatorname{nz}_0(N)^2+\operatorname{nz}_1(N)^2
\]
and
\[
\sum_{R=0}^N r(2R+1)\le 2\operatorname{nz}_0(N)\cdot\operatorname{nz}_1(N).
\]

Now we combine the three lines above with a variant of \eqref{eq:refinedlower} to get
\begin{align*}
2N+O(\operatorname{nz}(2N+1))&< \frac{1}{2}\left( \operatorname{nz}_0(N)^2+\operatorname{nz}_1(N)^2\right)+2\operatorname{nz}_0(N)\operatorname{nz}_1(N)\\
 &=\frac{1}{2}\operatorname{nz}(2N+1)^2-\operatorname{nz}_0(N)^2+\operatorname{nz}(2N+1)\cdot \operatorname{nz}_0(N).
\end{align*}
This is interesting because the last part of the inequality is maximized when $\operatorname{nz}_0(N)= \operatorname{nz}(2N+1)/2$, leading to
\[
2N+O(\operatorname{nz}(2N+1)) \le \frac{3}{4} \operatorname{nz}(2N+1).
\]
If the non-zero digits are not evenly distributed between even and odd indices, we could get even stronger results, although none of them would surpass the bound found in Theorem \ref{thm:main}. Perhaps a different subsequence or a clever weighting of the sum of $T(R)$'s would produce improved results.

\section{Comparing two different expansions}\label{sec:3root2}

The argument we have given for the number of non-zero digits in the expansion of $\sqrt{2}$ works just as well for bounding the number of non-zero digits in $3\sqrt{2}$. However, if we see the string $0100$ starting in the $n$th position in the expansion of $\sqrt{2}$, then we must see $1$'s at the $n$th and $n+1$st position in the expansion of $3\sqrt{2}$. 

So if both $\sqrt{2}$ and $3\sqrt{2}$ have close to the same number of $1$'s in their expansions, it must be because we see the strings $11$ or $101$ (or, possibly, $111$) appear in the expansion of $\sqrt{2}$ a lot. This is quite useful, since if $a_n, a_{n+k}, a_m, a_{m+k}=1$, with $n\neq m$, then a better than trivial bound can be placed on $r(n+m+k)$. In particular, it will count the pairs $(n,m+k),(m+k,n),(n+k,m),(m,n+k)$ and thus be at least $4$, so that $T(m+n+k-1)$ is at least $3$ by \eqref{eq:T(R-1)bound}, rather than the $2$ we typically assume. If one could show this happens often enough, one would get a non-trivial improvement in the lower bound.

However, in our attempts to use this technique, we could not do better than the results given in Theorem \ref{thm:2}, and so we leave it here as an idea in the hope that it inspires someone else to push the results further.

\end{document}